\documentclass[10pt]{amsart}
\usepackage{amsfonts}
\usepackage{graphicx}
\usepackage{amscd}
\usepackage{amsmath,amsfonts,amssymb,amsthm,mathrsfs,xypic}
\xyoption{curve}
\usepackage{fancybox}
\newcommand{\m}{\Lambda}

\newcommand{\cok}{\operatorname{Coker}}
\newcommand{\Ima}{\operatorname{Im}}

\newcommand{\ha}{\operatorname{\mathcal{A}}}

\newcommand {\hmc}{{\rm Mon}(Q, \mathcal{A} )}
\newcommand{\s}{\hfill \blacksquare}
\newtheorem{thm}{Theorem}[section]

\newtheorem{cor}[thm]{Corollary}

\newtheorem{lem}[thm]{Lemma}
\newtheorem{exm}[thm]{Example}

\newtheorem{prop}[thm]{Proposition}
\newtheorem{rem}[thm]{Remark}
\newtheorem{defn}[thm]{Definition}
\begin{document}

\title [Injective Objects of Monomorphism Categories ] {Injective Objects of Monomorphism Categories}
\author [Keyan Song, Zhanping Wang, Yuehui Zhang] {Keyan Song$^{1\dag}$, Zhanping Wang$^{1,2\dag\dag}$, Yuehui Zhang$^{1*}$}
\thanks{$^*$The corresponding author \ \ \ \ \ \ zyh@sjtu.edu.cn}
\thanks{Supported by the NSF of China (11271257, 11201377) and NSF of Shanghai (13ZR1422500).}
\thanks {$^\dag$ skyysh@outlook.com, $^{\dag\dag}$wangzp@nwnu.edu.cn}

\maketitle

\begin{center}
1. Department of Mathematics, \ \ Shanghai Jiao Tong University\\
Shanghai 200240, PR China\\
2. Department of Mathematics, \ \ Northwest Normal University\\
Lanzhou 730070, PR China
\end{center}

\begin{abstract}
\ For an acyclic quiver $Q$ and a finite-dimensional algebra $A$, we give a unified form of the indecomposable injective objects in the monomorphism category ${\rm Mon}(Q,A)$  and prove that ${\rm Mon}(Q, A)$ has enough injective objects. As applications, we show that for a given self-injective algebra $A$, a tilting object in the stable category $\underline{A}$-mod induces a natural tilting object in the stable monomorphism category $\underline{\rm Mon}(Q,A)$. We also realize the singularity category of the algebra $kQ\otimes_k A$ as the stable monomorphism category of the module category of $A$.

\vskip15pt


\vskip5pt

\end{abstract}

\section {\bf Introduction}

 For a $k$-algebra $A$ and a quiver $Q$ , the monomorphism category ${\rm Mon}(Q, A)$ is defined in [LZ]. It is a full subcategory of category ${\rm Rep}(Q,A)$ consisting of all the representations of $Q$ over $A$ (see section 2.2). Monomorphism category is interested by many people(cf. [B],[C1],[C2],[KlM1],[KLM2],\ [LZ],[Mo],[RS1], [RS2],[RS3],[SKZ],[XZZ],[Z]).

 When it comes to research this kind of category, a basic but important problem is to determine its projective objects and injective objects. Much research have been done about this problem:
 For $Q = \bullet \rightarrow \bullet$, all indecomposable projective (resp. injective) objects of ${\rm Mon}(Q, A)$ are given in [RS2]; For $Q = n \bullet \rightarrow \cdots \rightarrow
 \bullet 1$, all indecomposable projective (resp. injective) objects of ${\rm Mon}(Q, A)$ are determined by [Z] and [XZZ]. Recently, for any finite acyclic quiver $Q$, [LZ] has characterized all the indecomposable projective objects and one can easily see that ${\rm Mon}(Q, A)$ has enough projective objects.  However, up till now, there is not a paper to expound it in an overall way and one still does not know the concrete form of the injective objects of ${\rm Mon}(Q, A)$, let alone whether ${\rm Mon}(Q, A)$ has enough injective objects or not. This paper is mainly devoted to this problem (our interest still focus on finite acyclic quivers).

 As applications, we have three generalizations of the main results of [C1]. First of all, we prove that the monomorphism category ${\rm Mon}(Q, A)$ has enough injective objects for a given self-injective algebra $A$ (see Lemma \ref{enoughinj}). Secondly, we prove that for a given self-injective algebra $A$, a tilting object in the stable category $\underline{A}$-mod induces a natural tilting object in the stable monomorphism category $\underline{\rm Mon}(Q,A)$ (see Theorem \ref{tobj}). Finally, we relate the stable monomorphism category of the module category of an algebra $A$ to the singularity category of the algebra $kQ\otimes_k A$ (see Theorem \ref{sig}).

\vskip10pt

\section {\bf Monomorphism categories}

We fix notations and give necessary definitions and lemmas in this section.

\vskip10pt

\subsection{} Throughout this paper, $k$ is a field, $Q$ is a finite acyclic quiver (i.e., a finite quiver without oriented cycles), and
$A$ is a finite-dimensional $k$-algebra. Denote by $kQ$ the path
algebra of $Q$ over $k$. Let $P(i)$ (resp. $I(i)$) be the indecomposable projective
(resp. injective) $kQ$-module at $i\in Q_0$.  We denote the category of finite-dimensional
left $A$-modules by $A$-mod.

\vskip10pt

\subsection{} Given a finite acyclic quiver $Q = (Q_0,
Q_1, s, e)$, where $Q_0=\{1,2,\cdots, n\}$ is the set of vertices and $Q_1$ the set of
arrows. For any arrow $\alpha\in Q_1$, let $s(\alpha)$ and $e(\alpha)$ be the starting and the
ending point of $\alpha$ respectively.  Put $S$ the set of all source vertices (i.e.,  there is no arrow with $i$ as the ending point ) of $Q$. {\bf Note that we always assume that there is no arrow from $i$ to $j$ if $i< j$} in the sequel.

By definition (see [LZ]), {\it a
representation $X$ of $Q$ over $A$} is given by the following datum $X = (X_i, \
X_{\alpha}, \ i\in Q_0, \ \alpha \in Q_1)$, or simply $X = (X_i,
X_\alpha)$:

$\bullet$  each $X_i$ is an $A$-module;

$\bullet$ $X_{\alpha}:X_{s(\alpha)} \rightarrow X_{e(\alpha)}$ is an $A$-map.

 $X = (X_i,X_\alpha)$ is {\it a finite-dimensional representation} if each $X_i$ is finite-dimensional. A morphism $f$ from $X$ to $Y$
is a datum $(f_i, \ i\in Q_0)$, where $f_i: X_i \rightarrow Y_i$ is
an $A$-map for $i\in Q_0$, such that for each arrow $\alpha:
j\rightarrow i$ the following diagram

\[\xymatrix {X_j\ar[r]^{f_j}\ar[d]_{X_\alpha} & Y_j\ar[d]^{Y_\alpha} \\
X_i\ar[r]^{f_i} & Y_i}\] commutes. We call $f_i$ {\it the
$i$-th branch} of $f$. If $p = \alpha_l\cdots\alpha_1$
with $\alpha_i\in Q_1$,  $l\ge 1$, and $e(\alpha_i) =
s(\alpha_{i+1})$ for $1\le i\le l-1$, then we put $X_p$ to be the
$A$-map $X_{\alpha_l}\cdots X_{\alpha_1}$. Denote by ${\rm Rep}(Q,
A)$ the category of finite-dimensional representations of $Q$ over
$A$. Note that a sequence of representations $0\longrightarrow X\stackrel
{f}\longrightarrow Y\stackrel {g}\longrightarrow Z \longrightarrow
0$ in ${\rm Rep}(Q, A)$ is exact if and only if $0\longrightarrow
X_i\stackrel {f_i}\longrightarrow Y_i\stackrel {g_i}\longrightarrow
Z_i \longrightarrow 0$ is exact in $A$-mod for each $i\in Q_0$.

\vskip10pt

\begin{lem} \label{rep} {\rm ([LZ, Lemma 2.1])} \ We have an equivalence $\Lambda\mbox{-}{\rm
mod}\cong {\rm Rep}(Q, A)$ of categories.\end{lem}

In the following we will always identify a representation of $Q$ over $A$ with a $\Lambda$-module. If $T$ is an $A$-module and $M$ is a $kQ$-module
with $M = (M_i, i\in Q_0, M_\alpha, \alpha\in Q_1)\in {\rm Rep}(Q,
k)$, then $M\otimes_k T\in \m$-mod with $M\otimes_k T =
(M_i\otimes_kT = T^{{\rm dim}_kM_i}, \ i\in Q_0, \
M_\alpha\otimes_k{\rm Id}_T, \ \alpha\in Q_1)\in {\rm Rep}(Q, A)$(for more details, we refer [SKZ]).

\vskip10pt
 By ${\bf m}_i$ we denote the functor $P(i)\otimes_k - : \
 A{\rm -mod} \ \rightarrow {\rm Mon}(Q,A)$, and by ${\bf m}$ we
denote the functor $kQ\otimes_k-: A{\rm -mod}  \ \rightarrow {\rm Mon}(Q,A)$, where $P(i)$ is the $i$-th projective representation of quiver $Q$. Then ${\bf m}(M) = \bigoplus\limits_{\begin
{smallmatrix} i\in Q_0
\end{smallmatrix}} {\bf m}_i(M) = kQ\otimes_k M, \ \forall \ M\in A$-mod. It is clear that each functor ${\bf m}_i$ is exact and fully faithful.

We have the following useful lemma.
\begin{lem} \label{adj}  Let   $X=(X_j, X_\alpha) \in {\rm Rep}(Q, A)$ and $M \in A$-mod. Then we have an isomorphism of
abelian groups for each $i\in Q_0$, which is natural in both
positions
$${\rm Hom}_{{\rm Rep}(Q, A)}({\bf m}_i(M), X)\cong {\rm Hom}_{A}(M, X_i).$$
\end{lem}

 \subsection{} Here we recall from [LZ] the central notion of this paper.
\begin{defn} \label{maindef} \
{\rm ([LZ])}\ A representation $X = (X_i, X_{\alpha},
i\in Q_0, \alpha\in Q_1)\in {\rm Rep}(Q, A)$ is a monic
representation of $Q$ over $A$, or a monic $\Lambda$-module, if
$\delta_i(X)$ is an injective $A$-map for each $i\in Q_0$, where
$$\delta_i(X)=(X_{\alpha})_{\alpha\in Q_1, \ e(\alpha) = i}: \ \bigoplus\limits_{\begin
{smallmatrix} \alpha\in Q_1\\ e(\alpha) = i \end{smallmatrix}}
X_{s(\alpha)} \longrightarrow X_i.$$

Denote by ${\rm Mon}(Q, A)$ the full subcategory of ${\rm Rep}(Q,
A)$ consisting of all the monic representations of $Q$ over $A$,
which is called  the monomorphism category of $A$ over $Q$.
\end{defn}

\vskip5pt
\begin{rem}For $M\in kQ$-mod and $T\in A$-mod, it is clear that if
$M\in {\rm Mon}(Q, k)$ then $M\otimes_kT\in {\rm Mon}(Q, A)$. In
particular, $P(i)\otimes_k T\in {\rm Mon}(Q, A)$ for each $i\in
Q_0$. And that is why it is  easy for one to prove ${\rm Mon}(Q, A)$ has enough projective objects,
but it seems there is no such proof for one to show that ${\rm Mon}(Q, A)$ has enough injective objects (see [LZ], Proposition 2.4(i)).
\end{rem}
\vskip10pt

The following fact seems to be obvious, but the proof is quite involved in fact.

\begin{lem}\label{injobj}{\rm ([LZ, Proposition 2.4(iii)])} If $I$ is an indecomposable injective $A$-module and  $P(i)$
is the indecomposable projective $kQ$-module at $i\in Q_0$, then
$P(i)\otimes_k I$ is an indecomposable injective object in ${\rm
Mon}(Q, A)$.
\end{lem}

We also need the following obvious fact.

\begin{lem} \label{cok} Let $\mathcal{A}$ be an abelian category. If $f_i:X_i\rightarrow Y_i(i=1,2,\cdots,n)$ are morphisms in $\mathcal{A}$, then $\cok{\rm diag}\{f_i\}\cong \bigoplus_{i=1}^{n} \cok f_i$.
\end{lem}

\subsection{} ([SKZ]) Given $X = (X_j, X_\alpha)\in \m$-mod, for each $i\in Q_0$
we write $\cok\delta_i(X)$ (cf. Definition $\ref{maindef} $) as $\cok_i(X)$. Then we have a functor $\cok_i:
\m\mbox{-}{\rm mod}\longrightarrow A\mbox{-}{\rm mod},$
explicitly given by $\cok_i(X): = X_i/\sum\limits_{\begin {smallmatrix} \alpha\in Q_1\\
e(\alpha) = i \end{smallmatrix}} \Ima X_{\alpha}$ (if $i$ is a
source then $\cok_i(X): = X_i$). So we have an
exact sequence of $A$-modules $$\bigoplus\limits_{\begin
{smallmatrix} \alpha\in Q_1\\ e(\alpha) = i \end{smallmatrix}}
X_{s(\alpha)}\stackrel {\delta_i(X)}\longrightarrow  X_i \stackrel
{\pi_i(X)}  \longrightarrow \cok_i(X) \longrightarrow 0$$ for each
$X\in \m$-mod, where $\pi_i(X)$ is the canonical map.
\vskip10pt

\section{\bf Examples and Main Result}
We give three examples first, which will be useful for understanding Theorem \ref{mainthm}.
\begin{exm}\label{exm1} Let $Q$ be the quiver $2\rightarrow 1$. Then ${\rm Mon}(Q,\ A)$ has enough injective objects.
\end{exm}
\noindent{\bf Proof.} Let $X=X_2\stackrel{\delta_1}\rightarrow X_1$ be a monic representation in ${\rm Mon}(Q,\ A)$. Suppose that $f_2:X_2\hookrightarrow E_2$ and $\eta_1:\cok_1(X)\hookrightarrow K_1$ are injective envelopes of $X_2$ and $\cok_1(X)$ respectively. There exists some $s:X_1\rightarrow E_2$ such that $s\delta_1 =f_2$ since $\delta_1$ is monic and $E_2$ is injective. Put $f_1=\left(\begin{smallmatrix}s\\\eta_1\pi_1\end{smallmatrix}\right)$, $E_1=E_2\oplus K_1$ and $E=E_2\stackrel{(\begin{smallmatrix}1 \\ 0\end{smallmatrix})}\rightarrow E_1=(P(2)\otimes_k E_2)\oplus(P(1)\otimes_k K_1)$. It follows by Lemma \ref{injobj}, $E$ is injective. Consider the following commutative diagrams with exact rows
\[\xymatrix {
0\ar[r]&X_2\ar[d]^{f_2}\ar[r]^{\delta_1} &X_1\ar[d]^{f_1}\ar[r]^{\pi_1}&\cok_1(X)\ar[r]\ar[d]^{\eta_1}&0\\
0\ar[r]&E_2\ar[r]_{(\begin{smallmatrix}1 \\ 0\end{smallmatrix})} &E_2\oplus K_1\ar[r]_{(\begin{smallmatrix} 0\ 1\end{smallmatrix})}&K_1\ar[r]&0}
\]

By Snake lemma, $f_1$ is injective, and $\cok(f_1,\ f_2):X\rightarrow E$ is a monic representation.$\s$
\vskip5pt
\begin{rem} The above case is the lemma 2.1 of [C1].
\end{rem}

\begin{exm}\label{exm2} Let $Q$ be the quiver $1\leftarrow 3\rightarrow 2$. Then ${\rm Mon}(Q,\ A)$ has enough injective objects.
\end{exm}
\noindent{\bf Proof.}  Let $X=X_1\stackrel{\delta_1}\longleftarrow X_3\stackrel{\delta_2}\longrightarrow X_2$ be a monic representation in ${\rm Mon}(Q,\ A)$. Suppose that $f_3:X_3\hookrightarrow E_3$, $\eta_1:\cok_1(X)\hookrightarrow K_1$ and $\eta_2:\cok_2(X)\hookrightarrow K_2$ are injective envelopes of $X_2$, $\cok_1(X)$ and $\cok_2(X)$ respectively. It follows by Example \ref{exm1} that there exist two $A$-maps $f_1$, $f_2$ and monic representation $E=E_1\stackrel{(\begin{smallmatrix}1 \\ 0\end{smallmatrix})}\longleftarrow E_3\stackrel{(\begin{smallmatrix}1 \\ 0\end{smallmatrix})}\longrightarrow E_2$ such that $(f_1,\ f_2, \ f_3)$ is an injective morphism from $X$ to $E$, where $E_1=E_3\oplus K_1$ and $E_2=E_3\oplus K_2$. Obviously, $E=(P(3)\otimes_k E_3)\oplus(P(1)\otimes_k K_1)\oplus(P(2)\otimes_k K_2) $ is injective by Lemma \ref{injobj}. And by Snake lemma, $\cok(f_1,\ f_2, \ f_3)$ is a monic representation.$\s$

\begin{exm} \label{exm3}Let $Q$ be the quiver $2\rightarrow 1\leftarrow 3$. Then ${\rm Mon}(Q,\ A)$ has enough injective objects.
\end{exm}
\noindent{\bf Proof.}  Let $X=X_2\stackrel{\alpha}\longrightarrow X_1\stackrel{\beta}\longleftarrow X_3$ be a monic representation in ${\rm Mon}(Q,\ A)$. Suppose that $f_2:X_2\hookrightarrow E_2$, $f_3:X_3\hookrightarrow E_3$ and $\eta_1:\cok_1(X)\hookrightarrow K_1$ are injective envelopes of $X_2$, $X_3$ and $\cok_1(X)$ respectively. There is some $(\begin{smallmatrix} s_2 \\ s_3\end{smallmatrix}):X_1\rightarrow E_2\oplus E_3$ such that
$(\begin{smallmatrix} s_2 \\ s_3\end{smallmatrix})\delta_1=(\begin{smallmatrix} s_2 \\ s_3\end{smallmatrix})(\alpha\ \beta)=(\begin{smallmatrix} f_2\ 0 \\ 0\ f_3\end{smallmatrix})$ since $\delta_1$ is monic and $E_2\oplus E_3$ is injective. Put $f_1=\left(\begin{smallmatrix}s_2\\s_3\\\eta_1\pi_1\end{smallmatrix}\right)$, $E_1=E_2\oplus E_3\oplus K_1$ and $E=E_2\stackrel{(\begin{smallmatrix}1 \\ 0\\ 0\end{smallmatrix})}\longrightarrow E_1\stackrel{(\begin{smallmatrix}0 \\ 1\\ 0\end{smallmatrix})}\longleftarrow E_3$. It is clear that $E$ is a monic representation and an injective object in ${\rm Mon}(Q, A)$ since $ E=(P(1)\otimes_k K_1)\oplus(P(2)\otimes_k E_2)\oplus(P(3)\otimes_k E_3)$. Consider the following commutative diagrams with exact rows
\[\xymatrix {
0\ar[r]&X_2\oplus X_3\ar[d]^{(\begin{smallmatrix} f_2\ 0 \\ 0\ f_3\end{smallmatrix})}\ar[r]^{\delta_1} &X_1\ar[d]^{f_1}\ar[r]^{\pi_1}&\cok_1(X)\ar[r]\ar[d]^{\eta_1}&0\\
0\ar[r]&E_2\oplus E_3\ar[r]_{(\begin{smallmatrix}1\ 0\\ 0\ 1 \\ 0\ 0\end{smallmatrix})} &E_2\oplus E_3\oplus K_1\ar[r]_{(\begin{smallmatrix} 0 \ 0\ 1\end{smallmatrix})}&K_1\ar[r]&0}
\]

By Snake lemma, $f_1$ is injective, and $\cok(f_1,\ f_2, \ f_3):X\rightarrow E$ is a monic representation since $\cok(f_2)\oplus \cok(f_3)\cong \cok (\begin{smallmatrix} f_2\ 0 \\ 0\ f_3\end{smallmatrix})$ by Lemma \ref{cok}.$\s$

\vskip10pt
Inspired by Example \ref{exm1}, \ref{exm2} and \ref{exm3} we can prove our main result.

\begin{thm} \label{mainthm}Let $Q$ be a finite acyclic quiver and $A$ a finite-dimensional $k$-algebra. Then ${\rm Mon}(Q,\ A)$ has enough injective objects.
\end{thm}
\noindent{\bf Proof.} Let $X=(X_i,\ X_\alpha,\ i\in Q_0,\ \alpha\in Q_1)$ be a monic representation in ${\rm Mon}(Q,\ A)$. We divide our proof into two steps.
\vskip5pt
{\bf Step 1}. Construct injective monic representation $E=(E_i,\ E_\alpha,\ i\in Q_0,\ \alpha\in Q_1)$ and injective morphism $(f_1,\ f_2,\cdots, f_n):X\rightarrow E$.
\vskip5pt
{\it Substep 1:} Construct $E_i$ and $f_i:X_i\rightarrow E_i$ for each source $i\in S$.

  For each $i\in S$, let $f_i:X_i\rightarrow E_i$ be an injective envelope of $X_i$.
 \vskip5pt
  {\it Substep 2:} Construct $E_j$ and $f_j:X_j\rightarrow E_j$ for each $j\in Q_0-S$

 For each $j\in Q_0-S$,  we may assume that we have constructed $E_k$ and $f_k$ for any $k\in \{s(\alpha)|\alpha\in Q_1,\ e(\alpha)=j\}$. Recall that $$\delta_j(X)=(X_{\alpha})_{\alpha\in Q_1, \ e(\alpha) = j}: \ \bigoplus\limits_{\begin
{smallmatrix} \alpha\in Q_1\\ e(\alpha) = j \end{smallmatrix}}
X_{s(\alpha)} \longrightarrow X_j$$
and $\cok_j(X)=\cok_j(\delta_j(X))$. Let $\eta_j:\cok_j(X)\rightarrow K_j$ be an injective envelope of $\cok_j(X)$. Put $E_j= (\bigoplus\limits_{\begin
{smallmatrix} \alpha\in Q_1\\ e(\alpha) = j \end{smallmatrix}}
E_{s(\alpha)})\bigoplus K_j$. There is some $(t_{s(\alpha)})_{e(\alpha)=j}^\prime: X_j\rightarrow \bigoplus\limits_{\begin
{smallmatrix} \alpha\in Q_1\\ e(\alpha) = j \end{smallmatrix}}
E_{s(\alpha)}$ such that $$(t_{s(\alpha)})_{e(\alpha)=j}^\prime\delta_j(X)=(t_{s(\alpha)})_{e(\alpha)=j}^\prime(X_{\alpha})_{ \ e(\alpha) = j}={\rm diag}\{f_{s(\alpha)}\}_{e(\alpha)=j}$$ since $\delta_j(X)$ is monic and $\bigoplus\limits_{\begin
{smallmatrix} \alpha\in Q_1\\ e(\alpha) = j \end{smallmatrix}}
E_{s(\alpha)}$ is injective.

Let $f_j=((t_{s(\alpha)})_{e(\alpha)=j},\ \eta_j\pi_j)^\prime:X_j\rightarrow E_j$, and $\delta_j(E)$ be the canonical embedding. We have the following commutative diagrams with exact rows.
\[\xymatrix {
0\ar[r]&\bigoplus\limits_{\begin
{smallmatrix} \alpha\in Q_1\\ e(\alpha) = j \end{smallmatrix}}
X_{s(\alpha)}\ar[d]_{{\rm diag}\{f_{s(\alpha)}\}_{e(\alpha)=j}}\ar[r]^{\delta_j(X)} &X_j\ar[d]^{f_j}\ar[r]^{\pi_j(X)}&\cok_j(X)\ar[r]\ar[d]^{\eta_j}&0\\
0\ar[r]&\bigoplus\limits_{\begin
{smallmatrix} \alpha\in Q_1\\ e(\alpha) = j \end{smallmatrix}}
E_{s(\alpha)}\ar[r]^{\delta_j(E)} &E_j\ar[r]^{\pi_j(E)}&K_j\ar[r]&0}\eqno (3.1)
\]
By Snake lemma, $f_j$ is injective.

\vskip5pt

{\it Substep 3:} Construct $E_\alpha$ for each $\alpha\in Q_1$ and show $E=(E_i,\ E_\alpha,\ i\in Q_0,\ \alpha\in Q_1)$ is a monic representation.

For any $\alpha:i\rightarrow j\in Q_1$, $E_\alpha:E_i\rightarrow E_j$ is defined by $E_\alpha=\delta_j(E)\mu_i$, where $\mu_i:E_i\rightarrow E_j$ is the canonical embedding. Note that $\delta_j(E)$ is also the canonical embedding. By definition it is clear that $E=(E_i,\ E_\alpha,\ i\in Q_0,\ \alpha\in Q_1)\in {\rm Mon}(Q,A)$.

 \vskip5pt
  {\it Substep 4:} $E=(E_i,\ E_\alpha,\ i\in Q_0,\ \alpha\in Q_1)$ is an injective object in ${\rm Mon}(Q,A)$.

  Note that $E=(\bigoplus\limits_{i\in S} P(i)\otimes_k E_i)\bigoplus (\bigoplus\limits_{j\in Q_0-S} P(j)\otimes_k K_j))$, by Lemma \ref{injobj}, we see that $E$ is an injective object in ${\rm Mon}(Q,A)$.

\vskip5pt
{\it Substep 5:} $(f_1,\ f_2,\cdots, f_n):X\rightarrow E$ is an injective morphism in ${\rm Mon}(Q,A)$.

In fact, by the construction of {\it substep 1} and {\it substep 2} , we see that $(f_1,\ f_2,\cdots, f_n)$ is injective. We only need to prove that it is a morphism in ${\rm Mon}(Q,A)$. To this end, for each $\alpha:i\rightarrow j$, we need to prove the following diagram
 \[\xymatrix@R=0.5cm@C=0.8cm
{ X_i\ar[r]^-{X_\alpha}\ar[d]_{f_i} & X_j\ar[d]^{f_j}\\
E_i\ar[r]_-{E_\alpha} &\ E_j}\eqno (3.2)
\]
 commutes.

For convenience, put $j^-=\{s(\beta)|\ \beta \in Q_1, \ e(\beta)=j\}=\{i_1,\ i_2, \cdots, \ i_r\}$ and $\beta_{i_k}=i_k\rightarrow j$. Without losing of generality we can assume that $i=i_1$. So $X_{\beta_{i_1}}=X_\alpha$ and $E_{\beta_{i_1}}=E_\alpha$. Recall from {\it substep 2} that \begin{eqnarray*}(t_{i_1},t_{i_2}, \cdots, t_{i_r})^\prime \delta_j(X)=(t_{i_1},t_{i_2}, \cdots, t_{i_r})^\prime (X_{\beta_{i_1}},X_{\beta_{i_1}}, \cdots, X_{\beta_{i_1}})={\rm diag}\{f_{i_1},f_{i_2}, \cdots, f_{i_r}\}.\end{eqnarray*} It follows that

\begin{eqnarray*}f_jX_\alpha
&=&(t_{i_1},t_{i_2}, \cdots, t_{i_r},\ \eta_j\pi_j)^\prime X_{\beta_{i_1}}\\
&=&(t_{i_1}X_{\beta_{i_1}},t_{i_2}X_{\beta_{i_1}}, \cdots, t_{i_r}X_{\beta_{i_1}},\ \eta_j\pi_jX_{\beta_{i_1}})^\prime \\
&=&(f_{i_1},0,\cdots,0,0)^\prime\\
&=&E_{\beta_{i_1}}f_{i_1}\\
&=&E_\alpha f_{i}
\end{eqnarray*}
\vskip10pt
{\bf Step 2:} For above $(f_1, f_2, \cdots, f_n):X\rightarrow E$, $Y:=\cok(f_1, f_2, \cdots, f_n)$ is a monic representation.

 For each $i\in Q_0$, it is clear that $Y_i=\cok f_i$. For each $\alpha:i\rightarrow j\in Q_0$, $Y_\alpha:Y_i\rightarrow Y_j$ is induced by (3.2). Note that for any $j\in Q_0$, by Lemma \ref{cok} we have
\begin{eqnarray*}\cok {\rm diag}\{f_{s(\alpha)}\}_{e(\alpha)=j}=\bigoplus\limits_{\begin
{smallmatrix} \alpha\in Q_1\\ e(\alpha) = j \end{smallmatrix}}
\cok (f_{s(\alpha)})=\bigoplus\limits_{\begin
{smallmatrix} \alpha\in Q_1\\ e(\alpha) = j \end{smallmatrix}}
Y_{s(\alpha)}.
\end{eqnarray*}  Applying Snake lemma to diagram (3.1), we see $\delta_j(Y)=(Y_{\alpha})_{\alpha\in Q_1, \ e(\alpha) = j}: \ \bigoplus\limits_{\begin
{smallmatrix} \alpha\in Q_1\\ e(\alpha) = j \end{smallmatrix}}
Y_{s(\alpha)} \longrightarrow Y_j$ is injective.

This finishes the proof. $\s$

\vskip5pt

We have the following apparent corollary.

\begin{cor} The indecomposable injective objects in ${\rm Mon}(Q, A)$ are of the form $P(i)\otimes_k I$, where $P(i)$ is the indecomposable
projective $kQ$-module at $i\in Q_0$, and $I$ is an indecomposable injective $A$-module.
\end{cor}

\section{\bf Tilting Objects in Stable Monomorphism Category}
In this section, we assume that $A$ is a self-injective algebra. Thus $A$-mod is a Frobenius abelian category.
\subsection{}

As the first application of our main result, we have the following Lemma.
\begin{lem}\label{enoughinj} \ Let $A$ be a self-injective algebra. Then

${\rm (i)}$ \ ${\rm Mon}(Q, A)$
 is an exact category such that conflations are given by sequences $0\longrightarrow X\stackrel
{f}\longrightarrow Y\stackrel {g}\longrightarrow Z \longrightarrow
0$ with the induced sequence $0\longrightarrow
X_i\stackrel {f_i}\longrightarrow Y_i\stackrel {g_i}\longrightarrow
Z_i \longrightarrow 0$ exact in $A$-mod for each $i\in Q_0$.
\vskip5pt

${\rm (ii)}$ \ The exact category ${\rm Mon}(Q, A)$ is Frobenius such that its indecomposable projective objects are of the form ${\bf m}_i(P)$ where $P$ is an indecomposable projective module in  $A$-mod.
\end{lem}
\noindent{\bf Proof.} By Snake lemma, we see that ${\rm Mon}(Q, A)$ is closed under extension. Then it is clear that the exact sequences in ${\rm Mon}(Q, A)$ are induced by exact sequences in ${\rm Rep}(Q, A)$.

  And by [LZ] and Theorem \ref{mainthm}, ${\rm Mon}(Q, A)$ has enough projective objects, and its indecomposable projective objects are of the form ${\bf m}_i(P)$ where $P$ is an indecomposable projective  module in  $A$-mod. $\s$

\vskip10pt
 For the Frobenius abelian category $A$-mod, we denote by $\underline{{\rm Mon}}(Q, A)$ the stable category of
${\rm Mon}(Q, A)$ modulo projective objects; it is a triangulated category. We will call it the
{\it stable monomorphism category} of $A$-mod (for more details about Frobenius category we refer [H1]).
\vskip5pt
We have the following obvious proposition.
\begin{prop} For all $i\in Q_0$, ${\bf m}_i$ are fully-faithful and send projective objects to projective objects. Then they induce fully-faithful triangle functors ${\bf m}_i:\underline{A}{\rm -mod} \rightarrow \underline{{\rm Mon}}(Q, A)$.
\end{prop}

\begin{rem} If one replaces $A$-mod by a Frobenius abelian category $\ha$, then one can easily define {\it a representation $X$ of $Q$ over $\ha$} which is a datum $X
=(X_i, \ X_{\alpha}, \ i\in Q_0, \ \alpha \in Q_1)$, where $X_i$ is
an object in $\ha$  for each  $i\in Q_0$, and $X_{\alpha}: X_{s(\alpha)}
\longrightarrow X_{e(\alpha)}$ is a morphism in $\ha$ for each $\alpha \in
Q_1$. The monic representation and monomorphism category $\hmc$ can also be defined. By this way, the lemma 2.1 of [C1] can be generalized.
\end{rem}
\vskip10pt

\subsection{}In this section, we will show that for a self-injective algebra $A$, a tilting object $T$
in the stable category $\underline{A}{\rm -mod} $ induces a natural tilting object $kQ\otimes_k T$ in the stable monomorphism
category $\underline{{\rm Mon}}(Q, A)$.

 Let $\mathcal{T} $ be an algebraical triangulated category(it is triangle equivalent to the stable category of a Frobenius exact category). Denote by $[1]$ the shift functor and by $[n]$ its $n$-th power for each $n\in \mathbb{Z}$. An object $T$ in $\mathcal{T} $  is a tilting object if the following conditions are satisfied([C1]):

(T1) ${\rm Hom}_\mathcal{T} (T, T[n]) = 0 $ for $n \neq 0$;

(T2) the smallest thick(closed under taking summands) triangulated subcategory of $\mathcal{T} $ containing $T$ is $\mathcal{T} $ itself;

(T3) ${\rm End}_\mathcal{T} (T)$ is an artin algebra with finite global dimension.

\vskip10pt

The following is our second main application.
\begin{thm}\label{tobj} Let $A$ be a self-injective algebra and $T$ a tilting object
in its stable category $\underline{A}{\rm -mod} $. Then $kQ\otimes_k T$ is a tilting object in $\underline{{\rm Mon}}(Q, A)$;
Moreover, we have an isomorphism ${\rm End}_{\underline{{\rm Mon}}(Q, A)} (kQ\otimes_k T)\simeq kQ\otimes_k {\rm End}_{\underline{A}{\rm -mod} }(T)$ of algebras.
\end{thm}

\subsection{} Before proving Theorem \ref{tobj}, we first give the following lemma.
\begin{lem}\label{tlem} For each object $X=(X_i, X_\alpha)\in {\rm Mon}(Q, A)$, there is an exact sequence
$$0\rightarrow \bigoplus\limits_{\begin
{smallmatrix} i\in S
\end{smallmatrix}} P(i)\otimes_k X_i\stackrel{\phi}\rightarrow X\rightarrow Z\rightarrow 0$$ such that $Z=\cok(\phi)\in {\rm Mon}(Q, A)$.
\end{lem}
\noindent{\bf Proof.} For convenience, put $Y=\bigoplus\limits_{\begin
{smallmatrix} i\in S
\end{smallmatrix}} P(i)\otimes_k X_i$ which has the following description:

$\bullet$ $Y_j=\bigoplus\limits_{\begin{smallmatrix} i\in S\end{smallmatrix}}X_i^{m_{ij}}$, where $m_{ij}$ is the number of paths from $i$ to $j$;

$\bullet$ For each $\alpha:j\rightarrow k \in Q_1$, $Y_\alpha:Y_j\rightarrow Y_k$ is the natural embedding.

For each $j\in Q_0$, put $\delta_j=(X_p)_{s(p)=i, e(p)=j}$. It is clear that $\delta_j$ is injective. Then we have the following exact sequence
$$0\rightarrow \bigoplus\limits_{\begin{smallmatrix} i\in S\end{smallmatrix}}X_i^{m_{ij}}\stackrel{\delta_j}\rightarrow X_j\stackrel{\pi_j}\rightarrow Z_j\rightarrow 0$$
where $Z_j=\cok (X_p)_{s(p)=i,e(p)=j}$ and $\pi_j$ is canonical.

Obviously, $Y_j=\bigoplus\limits_{\begin{smallmatrix} i\in S\end{smallmatrix}}X_i^{m_{ij}}=X_j$ and $Z_j=\cok (X_p)_{s(p)=i,e(p)=j}=0$ provided $j\in Q_0$ is a source vertex.

For each $\alpha:j\rightarrow k\in Q_1$, we have the following diagram with exact rows

\[\xymatrix {
0\ar[r]&\bigoplus\limits_{\begin{smallmatrix} i\in S\end{smallmatrix}}X_i^{m_{ij}}\ar[d]^{Y_\alpha}\ar[r]^{\delta_j} &X_j\ar[d]^{X_\alpha}\ar[r]^{\pi_j}&Z_j\ar[r]\ar[d]^{Z_\alpha}&0\\
0\ar[r]&\bigoplus\limits_{\begin{smallmatrix} i\in S\end{smallmatrix}}X_i^{m_{ik}}\ar[r]_{\delta_k} &X_k\ar[r]_{\pi_j}&Z_k\ar[r]&0}\eqno (5.1)
\]

We need to show that $Z=(Z_i, Z_\alpha)\in {\rm Mon}(Q, A)$. For each $k\in Q_0$ (we may assume that $k$ is not a source vertex). Denote by $k^-=\{j_1,j_2,\cdots, j_s\}$ the set of all direct predecessors of $k$. Consider the following diagram with exact rows

\[\xymatrix {
0\ar[r]&\bigoplus\limits_{\begin{smallmatrix} r=1\end{smallmatrix}}^s\bigoplus\limits_{\begin{smallmatrix} i\in S\end{smallmatrix}}X_i^{m_{ij_s}}\ar[d]^{(Y_\alpha)_{s(\alpha)\in k^-}}\ar[r]_{\delta} &\bigoplus\limits_{\begin{smallmatrix} r=1\end{smallmatrix}}^s X_{j_r}\ar[d]^{(X_\alpha)_{s(\alpha)\in k^-}}\ar[r]_{\pi}&\bigoplus\limits_{\begin{smallmatrix} r=1\end{smallmatrix}}^s Z_{j_r}\ar[r]\ar[d]^{(Z_\alpha)_{s(\alpha)\in k^-}}&0\\
0\ar[r]&\bigoplus\limits_{\begin{smallmatrix} i\in S\end{smallmatrix}}X_i^{m_{ik}}\ar[r]_{\delta_k} &X_k\ar[r]_{\pi_k}&Z_k\ar[r]&0}\eqno (5.2)
\]
where $\delta={\rm diag}\{\delta_{j_1},\cdots,\delta_{j_s} \}$ and  $\pi={\rm diag}\{\pi_{j_1},\cdots,\pi_{j_s} \}$.

Note that $(Y_\alpha)_{s(\alpha)\in k^-}$ is actually the identity map. By Snake lemma, $(Z_\alpha)_{s(\alpha)\in k^-}$ is injective. It follows that $Z=(Z_i, Z_\alpha)$ is a monic representation.$\s$

\vskip10pt

If we delete the source vertices of a quiver $Q$, then we get a new quiver $Q^\prime$. Similarly, for any $X\in {\rm Rep}(Q, A)$, we obtain a representation $X^\prime \in {\rm Rep}(Q^\prime, A)$.

We have the following direct corollary.
\begin{cor}\label{tcor} Let $Q'$ be the quiver obtained by deleting the source vertices of a quiver $Q$ and $Z$ the one in \ref{tlem}. Then there is an exact sequence
$$0\rightarrow \bigoplus\limits_{\begin
{smallmatrix} i\in S(Q^\prime)
\end{smallmatrix}} P'(i)\otimes_k Z_i^\prime\stackrel{\phi^\prime}\rightarrow Z^\prime\rightarrow \cok({\phi^\prime})\rightarrow 0$$ such that $\cok(\phi^\prime)\in {\rm Mon}(Q^\prime, A)$, where $P'(i)$ is the projective $kQ'$-module at vertex $i$ and $S(Q^\prime)$ is the set of source vertices of quiver $Q^\prime$.
\end{cor}
\vskip10pt

\subsection{\bf Proof of Theorem 4.4}%\ref{tobj}}

Step 1. Check (T1) for $kQ\otimes_k T$.

 Note that $kQ\otimes_k T={\bf m}(T)=\bigoplus\limits_{\begin
{smallmatrix} i\in Q_0
\end{smallmatrix}} {\bf m}_i(T) =\bigoplus\limits_{\begin
{smallmatrix} i\in Q_0
\end{smallmatrix}} P(i)\otimes_k T$ and ${\bf m}_i(T[n])\simeq {\bf m}_i(T)[n]$. For any $j> i$, it follows by Lemma \ref{adj} that
$${\rm Hom}_{{\rm Rep}(Q, A)}({\bf m}_j(T), {\bf m}_i(T))\cong {\rm Hom}_{A}(T,({\bf m}_i(T))_j )={\rm Hom}_{A}(T,0 )=0.$$

So it suffices to show that for any $j\leq i$ and $n\neq 0$, we have $${\rm Hom}_{{\rm Rep}(Q, A)}({\bf m}_j(T), {\bf m}_i(T[n]))=0.$$

By definition, each $f\in {\rm Hom}_{\underline{A}}(T, T[n])$ factors through some projective object $P$ in $A$-mod, it follows that $\phi$ factors through ${\bf m}_j(P)$ for any $\phi\in {\rm Hom}_{{\rm Rep}(Q, A)}({\bf m}_j(T), {\bf m}_i(T[n]))$.

Step 2. Check (T2) for $kQ\otimes_k T$.

By Lemma \ref{tlem}, each object $X\in {\rm Mon}(Q, A)$ fits into a conflation
$$0\rightarrow \bigoplus\limits_{\begin
{smallmatrix} i\in S
\end{smallmatrix}} {\bf m}_i( X_i)\rightarrow X\rightarrow Z\rightarrow 0$$

\noindent and thus into a triangle $$ \bigoplus\limits_{\begin
{smallmatrix} i\in S
\end{smallmatrix}} {\bf m}_i (X_i)\rightarrow X\rightarrow Z\rightarrow (\bigoplus\limits_{\begin
{smallmatrix} i\in S
\end{smallmatrix}} {\bf m}_i (X_i))[1]$$

Using Corollary \ref{tcor} and Lemma \ref{tlem} alternatively, we see that the smallest triangulated
subcategory of $\underline{{\rm Mon}}(Q, A)$ containing all ${\bf m}_i(\underline{A}-mod)$ is $\underline{{\rm Mon}}(Q, A)$ itself. Now applying
the condition (T2) of $T$, we infer that (T2) holds for $kQ\otimes_k T$.

Step 3. Check (T3) for $kQ\otimes_k T$.

It is direct to check that ${\rm End}_{\underline{{\rm Mon}}(Q, A)} (kQ\otimes_k T)\simeq kQ\otimes_k {\rm End}_{\underline{A}-mod}(T)$.
 Recall that the algebra ${\rm End}_{\underline{A}-mod}(T)$ has finite global dimension. Then
by ([ARS], Chapter III, Proposition 2.6) we infer that ${\rm End}_{\underline{{\rm Mon}}(Q, A)} (kQ\otimes_k T)$ has finite global dimension.

\section{\bf Stable Monomorphism Category as Singularity Category}

In this section, we will relate the stable monomorphism category of the module category of a finite-dimensional $k$-algebra $A$ to the singularity category of the algebra $kQ\otimes_k A$.

\subsection{} We first give the following definition.

\begin{defn}\label{maindef'}{\rm [SKZ]} Let $\mathcal X$ be an  additive full  subcategory of
$A$-mod. Denote by ${\rm Mon}(Q, \mathcal X)$ the full subcategory
of ${\rm Mon}(Q, A)$ consisting of all the monic representations $X
= (X_i, X_{\alpha})$, such that $X_i\in \mathcal X$ and
$\cok\delta_i(X) \in \mathcal X$ for all $i\in Q_0$. We call ${\rm
Mon}(Q, \mathcal X)$ the monomorphism category of $\mathcal X$ over
$Q$.
\end{defn}
\subsection{} Modules in $^\perp A$ are called Cohen-Macaulay modules. Denote  $^\perp A=\mathbf{CM}(A)$.
  {\it A complete $A$-projective
resolution} is an exact sequence of finitely generated projective
$A$-modules
$$P^\bullet = \cdots \longrightarrow P^{-1}\longrightarrow P^{0}
\stackrel{d^0}{\longrightarrow} P^{1}\longrightarrow \cdots$$ such
that ${\rm Hom}_A(P^\bullet, A)$ remains to be exact. A module $M\in
A$-mod is {\it Gorenstein-projective}, if there exists a complete
$A$-projective resolution $P^\bullet$ such that $M\cong
\operatorname{Ker}d^0$. Let $A-\mathcal{G}proj$ be the full subcategory of $A$-mod of  Gorenstein-projective $A$-modules.
Then $A-\mathcal{G}proj\subseteq \mathbf{CM}(A)$; and if $A$ is Gorenstein algebra(i.e, the right and left injective dimension of $A$ are both finite), then  $A-\mathcal{G}proj= \mathbf{CM}(A)$([EJ], Corollary 11.5.3).

\begin{lem} \label{skz}{\rm [SKZ]} Let $A$ be a finite-dimensional k-algebra and $T$ an $A-$module. If there exists an exact sequence
\[\begin{CD} 0@>>>T_m@>>>\cdots @>>>T_0@>>> D(A_A)@>>>0\end{CD}\]
such that $T_i\in {\rm add }\ T$, then $^\perp (kQ\otimes_k T)={\rm Mon}(Q, ^\perp T).$
\end{lem}

\begin{cor} Let $A$ be a finite-dimensional k-algebra with inj.dim $A_A< \infty $.

{\rm (i)} Then $\mathbf{CM}(kQ\otimes_k A)=\rm{Mon}(Q,\mathbf{CM}(A)).$

{\rm (ii)} Assume further that $A$ is Gorenstein. Then $kQ\otimes_kA$-$\mathcal{G}proj$=$\rm{Mon}(Q, A-\mathcal{G}proj)$.
\end{cor}

\noindent{\bf Proof.} If inj.dim $A_A< \infty $, we can take $T=A$ in Lemma \ref{skz}, then we have $\rm{Mon}$$(Q,^\perp A)=^\perp (kQ\otimes_k A)$, so we proved (1). If $A$ is Gorenstein algebra, then it is well known that $kQ\otimes_k A$ is again Gorenstein, thus (2) follows from (1).

\subsection{} Let $\mathcal{D}^b(A)$ be the bounded derived category of $A$, and $\mathcal{K}^b(\mathcal{P}(A))$ be the bounded homotopy category of $\mathcal{P}(A)$. The singularity category $\mathcal{D}_{sg}^b(A)$ of $A$ is defined to be the Verdier quotient $\mathcal{D}^b(A)/\mathcal{K}^b(\mathcal{P}(A))$. If $A$ is Gorenstein, then there is a triangle-equivalence $\mathcal{D}_{sg}^b(A)\cong \underline{\mathbf{CM}(A)}$, where $\underline{\mathbf{CM}(A)}$ is the stable category of $\mathbf{CM}(A)$ ([H2],Theorem 4.6).

\vskip5pt
As our third application, it follows by Corollary 5.2 we have the following main theorem of this section.

\begin{thm}\label{sig}  Let $A$ be a Gorenstein algebra. Then there is a triangle-equivalence
 $$\underline{\rm{Mon}}(Q,\mathbf{CM}(A))\simeq \mathcal{D}_{sg}^b(kQ\otimes_kA).$$
 In particular, if $A$ is a self-injective algebra, then $$\underline{\rm{Mon}}(Q,A)\simeq \mathcal{D}_{sg}^b(kQ\otimes_kA).$$
\end{thm}

\end{document}